\documentclass{amsart}
\usepackage{graphicx}
\usepackage{mathtools}
\usepackage{amsmath}
\usepackage{bbm}
\newcommand{\indicator}{\mathbbm{1}}
\DeclareMathOperator*{\essinf}{ess\,inf}
\DeclareMathOperator*{\esssup}{ess\,sup}
\newtheorem{thm}{Theorem}[section]
\newtheorem{cor}[thm]{Corollary}

\theoremstyle{definition}
\newtheorem{defn}[thm]{Definition}
\theoremstyle{remark}

\numberwithin{equation}{section}

\begin{document}
	
	\title[Schauder bases in Analytic variable exponent Hardy spaces]{Schauder bases in Analytic variable exponent Hardy spaces}%
	
	\author{Tengiz Kopaliani}
	\address{Tengiz Kopaliani \\
		Faculty of Exact and Natural Sciences\\
		Javakhishvili Tbilisi State
		University\\
		13, University St., Tbilisi, 0143, Georgia}
	\email{tengiz.kopaliani@tsu.ge}
	
	\author{Giorgi Machitidze}
	\address{
		Giorgi Machitidze \\
		Faculty of Exact and Natural Sciences\\
		Javakhishvili Tbilisi State
		University\\
		13, University St., Tbilisi, 0143, Georgia}
	\email{ <giorgimachitidze13@gmail.com>}
	\address{}%
	\email{}%
	\thanks{}%
	
	
	\begin{abstract}
		We study Schauder bases in analytic variable exponent Hardy spaces with $1<p_{-}\leq p_{+}<\infty$. We prove that, for a sequence $(a_n)_{n\geq1}\subset\mathbb D$ satisfying $\sum_{n=1}^{\infty}(1-|a_n|)=\infty$, the corresponding Malmquist--Takenaka system is a Schauder basis in $H^{p(\cdot)}(\mathbb T)$ if and only if the exponent satisfies the variable Muckenhoupt condition $p(\cdot)\in\mathcal A(\mathbb T)$. Under the additional condition $\sup_{n\geq1}|a_n|<1$, the same criterion holds for the corresponding system of finite Blaschke products.
	\end{abstract}
	\maketitle
	
	\section{Introduction}
	
	Let $\mathbb{D}$ be the open unit disk of the complex plane. We denote by $\mathbb{T}=\{z\in \mathbb{C}:\,\,|z|=1\}$ the unit circle and by $dm$ the normalized Lebesgue measure on $\mathbb{T}$ (most often identified with the normalized Lebesgue measure $\frac{dx}{2\pi}$ on the
	interval $[0, 2\pi)$).
	
	It is well known that the classical $H^{p}$ $(1\leq p\leq\infty)$ spaces can be defined in
	two equivalent ways:
	
	1) $H^{p}(\mathbb{D})$ is the space of analytic functions $f:\mathbb{D}\rightarrow \mathbb{C}$ for which, setting $f_{r}(t)=f(re^{it})$,
	$$
	\|f\|_{H^{p}(\mathbb{D})}=\sup_{0\leq r<1}\Arrowvert f_{r}\Arrowvert_{p}
	$$
	is finite. The Fatou-Riesz Theorem asserts that the boundary limits
	$$
	f^{\ast}(t)=\lim_{r\rightarrow 1}f_{r}(t)
	$$
	exist almost everywhere and $\|f\|_{H^{p}(\mathbb{D})}=\|f^{\ast}\|_{p}$. One has $f^{\ast}\in L^{p}(\mathbb{T})$, and its Fourier coefficients $\widehat{f}(n)=0$ for $n<0.$
	
	2) Conversely, for every function $g\in H^{p}(\mathbb{T})$ ($g\in L^{p}(\mathbb{T})$ with Fourier coefficients $\widehat{g}(n)=0$ for $n<0$), the analytic extension $P[g]:\mathbb{D}\rightarrow \mathbb{C}$ defined by $P[g](z)=\sum_{n\geq 0}\widehat{g}(n)z^{n}$ is in $H^{p}(\mathbb{D})$ and $\Arrowvert P[g]\Arrowvert_{H^{p}(\mathbb{D})}=\Arrowvert g\|_{p}.$
	
	In most cases, by convention, we may identify a function $f\in H^{p}(\mathbb{D})$ with its boundary values $bf\in H^{p}(\mathbb{T}) $ and hence also the spaces $H^{p}(\mathbb{D})$ and $H^{p}(\mathbb{T})$.
	
	Analytic variable exponent Hardy spaces on the unit disk extend classical Hardy spaces by replacing a constant integrability exponent $p$ with a measurable function $p(\cdot)$. In the case when the exponent $p(\cdot)$ satisfies the log-H\"{o}lder condition, some properties of the space $H^{p(\cdot)}(\mathbb{D})$ were
	investigated in \cite{KP1, KP2} and \cite{CHCH}. In the present paper, we consider the $H^{p(\cdot)}(\mathbb{D})$ (and $H^{p(\cdot)}(\mathbb{T})$) spaces in a more general setting, namely when the exponent $p(\cdot)$ satisfies the Muckenhoupt condition for variable exponent spaces (in short, $p(\cdot)\in \mathcal{A}(\mathbb{T}) $). In the context of harmonic analysis tools, the variable analogue of the classical Muckenhoupt condition plays a pivotal role. This condition, in the case $1<p_{-}\leq p_{+}\leq\infty$, is a necessary and sufficient condition for the boundedness of the Hardy-Littlewood maximal operator in $L^{p(\cdot)}(\mathbb{T}).$ In this article, we introduce a version of the variable exponent Hardy spaces $H^{p(\cdot)}(\mathbb{D})$ and $H^{p(\cdot)}(\mathbb{T})$. We identify these spaces in the case $p(\cdot)\in \mathcal{A}(\mathbb{T})$. We consider problems related to the existence of bases of rational functions in the $H^{p(\cdot)}(\mathbb{D})$ space (indeed, bases of finite Blaschke products or the Malmquist-Takenaka system).
	
	Let $ X$ be a complex topological linear space. We say that the sequence $(x_{n})_{n=1}^{\infty}$ in $X$ is
	a Schauder basis for $X$ if, for each $x \in X$, there is a unique sequence $(a_{n})_{n=1}^{\infty}$ of complex numbers such that $x = \sum_{n=1}^{\infty} a_{n}x_{n}$, where the series converges in the topology of $X$.
	It is well known that the trigonometric system $\left\lbrace \frac{e^{ikt}}{\sqrt{2\pi}} \right\rbrace_{k=-\infty}^{\infty} $ is a Schauder basis in $L^{p}(\mathbb{T})$ for all $1<p<\infty$, and the system $\left\lbrace \frac{e^{ikt}}{\sqrt{2\pi}} \right\rbrace_{k=0}^{\infty} $ is a Schauder basis in the Hardy space $H^{p}(\mathbb{T})$ for all $1<p<\infty$.
	
	For a given system of points $(a_{n})_{n=1}^{\infty}$ inside the unit disk (the points $a_{k}$ are enumerated taking into account
	their multiplicity), the Takenaka–Malmquist system (TM system)
	is defined by
	$$
	\phi_{n}(z)=B_{n}(z)\frac{\sqrt{1-|a_{n+1}|^{2}}}{1-\overline{a_{n+1}}z},\,\,\,n\in \mathbb{Z}_{+},
	$$
	where
	$$
	B_{n}(z)=\prod_{j=1}^{n}\frac{z-a_{j}}{1-\overline{a_{j}}z},\,\,n\in \mathbb{N}\,\,\,\mbox{and}\,\,B_{0}=1
	$$
	is a Blaschke product for the disk $\mathbb{D}$ of degree $n$.
	
	It is well known that the TM system is an orthonormal and complete system in $H^{2}(\mathbb{T})$ if and only if $B_{n}(z)\rightarrow 0$ as $n\rightarrow\infty $
	uniformly on compact sets in $\mathbb{D}$. The last condition is equivalent to non-Blaschke condition (see \cite{Bur})
	\begin{equation}\label{Blaschke condition}
		\sum_{n=1}^{\infty}(1-|a_{n}|)=+\infty.
	\end{equation}

	For fixed non-Blaschke sequence in $\mathbb{D}$, we denote corresponding TM system  by $(\phi_n)_{n\geq0}$ and sequence of finite Blaschke products by  $(B)_{n\geq0}$ (B system).  
	If $a_{n}\equiv0,$ the corresponding TM system reduces to the classical trigonometric system. The Malmquist-Takenaka system is a perturbation of the classical trigonometric system, where powers of $z$ are replaced by products of other M\"{o}bius transforms of the disc. 
	Note that TM system is intimately connected with popular non-linear phase unwinding decomposition techniques in both mathematical analysis and signal processing (see \cite{CP, TO1, TO2}).  It was shown by Coifman and Peyri\'{e}re in \cite{CP} and by Qian, Chen, and Tan in \cite{TQL} that $(\phi_n)_{n\geq0}$ is a basis in the Hardy spaces	$H^{p}(\mathbb T)$, $1<p<\infty.$  
	 Recently, it was shown in \cite{FMN} that if the sequence \((a_n)_{n=1}^\infty\) is contained in a compact subset of the open unit disk, then \((B_n)_{n\geq 0}\) is a Schauder basis for \(\text{Hol}(\overline{\mathbb{D}})\) endowed with its natural topology. Furthermore, the authors also fully characterized all sequences \((a_n)_{n=1}^\infty\) for which the corresponding sequence of finite Blaschke products forms a Schauder basis in \(H^p(\mathbb{T})\) \((1 < p < \infty)\).

	The main results of the present paper are the following theorems.
	
	\begin{thm}\label{main Theorem1}
		Let $1 < p_{-} \leq p_{+} < \infty$ and $(a_n)_{n\geq1}\subset\mathbb D$
		satisfy \eqref{Blaschke condition}.
		Let $(\phi_n)_{n\geq0}$ be the corresponding TM
		system. Then the following assertions are equivalent:
		\begin{enumerate}
			\item[(i)] $p(\cdot)\in\mathcal A(\mathbb T)$;
			\item[(ii)] $(\phi_n)_{n\geq0}$ is a Schauder basis of
			$H^{p(\cdot)}(\mathbb T)$.
		\end{enumerate}
	\end{thm}
	
	\begin{thm}\label{equivalence}
		Let $1<p_{-}\leq p_{+}<\infty$ and let $(a_n)_{n\geq1}\subset\mathbb D$ satisfy
		\begin{equation}\label{uniform-away-boundary}
			\sup_{n\geq1}|a_n|<1.
		\end{equation}
		Let $(\phi_n)_{n\geq0}$ and $(B_n)_{n\geq0}$ be the corresponding TM and  B systems. Then the following assertions are equivalent:
		\begin{enumerate}
			\item[(i)] The system $(B_n)_{n\geq0}$ is a Schauder basis in $H^{p(\cdot)}(\mathbb T)$;
			\item[(ii)] $(\phi_n)_{n\geq0}$ is a Schauder basis in $H^{p(\cdot)}(\mathbb T).$
		\end{enumerate}
	\end{thm}
	
	Note that from Theorem \ref{main Theorem1} and Theorem \ref{equivalence} we obtain:
	
	\begin{thm}\label{main Theorem2}
		Let $1<p_{-}\leq p_{+}<\infty$ and let $(a_n)_{n\geq1}\subset\mathbb D$ satisfy
		\begin{equation}\label{uniform-away-boundary}
			\sup_{n\geq1}|a_n|<1.
		\end{equation}
		Then the following assertions are equivalent:
		\begin{enumerate}
			\item[(i)] $p(\cdot)\in\mathcal A(\mathbb T)$;
			\item[(ii)] the system $(B_n)_{n\geq0}$ is a Schauder basis of
			$H^{p(\cdot)}(\mathbb T)$.
		\end{enumerate}
	\end{thm}

	\section{Some classical operators in the $L^{p(\cdot)}(\mathbb{T})$ space}
	
	A measurable function $p: \mathbb{T} \to [1, \infty]$ is called a variable exponent on $\mathbb{T} $. We denote the set of all such variable exponents by $\mathcal{P}(\mathbb{T} )$. For a given $p \in \mathcal{P}(\mathbb{T})$, we define
	\begin{align*}
		p_{-} \coloneqq \essinf_{x \in \mathbb{T} } p(x) \qquad \text{and} \qquad p_{+} \coloneqq\esssup_{x \in \mathbb{T} } p(x).
	\end{align*}
	The conjugate exponent $p' \in \mathcal{P}(\mathbb{T} )$ is defined pointwise by the relation $\frac{1}{p(x)} + \frac{1}{p'(x)} = 1$.
	
	For a measurable function $f: \mathbb{T} \to \mathbb{C}$, we define the modular
	\begin{align*}
		\rho_{p(\cdot)}(f) \coloneqq \int_{\mathbb{T} } |f(x)|^{p(x)} \, dm(x),
	\end{align*}
	where we use the standard convention $t^\infty \coloneqq 0$ for $t \in \mathbb{T} $ and $t^\infty \coloneqq \infty$ for $t > 1$. The variable exponent Lebesgue space $L^{p(\cdot)}(\mathbb{T} )$ is defined as the set of all measurable functions $f$ such that $\rho_{p(\cdot)}(f/\lambda) < \infty$ for some $\lambda > 0$. It becomes a Banach space when equipped with the Luxemburg norm
	\begin{align*}
		\|f\|_{p(\cdot)} \coloneqq \inf \left\lbrace \lambda > 0 : \rho_{p(\cdot)}(f/\lambda) \leq 1 \right\rbrace.
	\end{align*}
	If the exponent $p(\cdot)$ is constant, $L^{p(\cdot)}(\mathbb{T} )$ simply coincides with the classical Lebesgue space. Moreover, the generalized Hölder inequality holds in these spaces: for any $f \in L^{p(\cdot)}(\mathbb{T} )$ and $g \in L^{p'(\cdot)}(\mathbb{T} )$, we have
	\begin{align*}
		\int_{\mathbb{T} } |f(x)g(x)| \, dmx \leq 2 \|f\|_{p(\cdot)} \|g\|_{p'(\cdot)}.
	\end{align*}
	
	The duality pairing can be written as
	$$
	<f,g>=\int_{\mathbb{T}}f(x)\overline{g(x)}dm(x)
	$$
	where $f\in L^{p(\cdot)}(\mathbb{T}) $ and $g\in L^{p'(\cdot)}(\mathbb{T}).$
	Given $f\in L^{1}(\mathbb{T})$, the Hardy-Littlewood maximal function is defined by
	$$
	Mf(t)=\sup_{t \in I}\frac{1}{m(I)}\int_{I}|f(x)|dm(x),\,\,\,t\in \mathbb{T},
	$$
	where the supremum is taken over all arcs $I \subset \mathbb{T}$ containing $t\in \mathbb{T}$.
	
	Let
	$$
	P_{r}(f,t)=\int_{\mathbb{T}}\dfrac{(1-r^{2})f(y)dm(y)}{1-2rcos(t-y)+r^{2}}\,\,\,t\in \mathbb{T}
	$$
	be the Poisson integral of $f$.
	
	The circular Hilbert transform is defined by
	$$
	Hf(t)=\operatorname{p.v.}\int_{\mathbb{T}}\frac{f(y)dm(y)}{\tan\frac{t-y}{2}} \,\,,t\in \mathbb{T}.
	$$
	
	Let $P_{+}$ be the M. Riesz projection of $L^{p}(\mathbb{T})$ onto $H^{p}(\mathbb{T})$ ($1<p<\infty$), defined by
	$$
	P_{+}\left( \sum_{-\infty}^{\infty}a_{n}e^{int}\right) =\sum_{0}^{\infty}a_{n}e^{int}
	$$
	where the sum is the Fourier series of an arbitrary $f\in L^{p}(\mathbb{T})$.
	
	In the context of harmonic analysis tools, the variable analogue of the classical Muckenhoupt condition plays a pivotal role. For a variable exponent $p \in \mathcal{P}(\mathbb{T} )$, we say that $p$ satisfies the $\mathcal{A}$ condition (and write $p(\cdot) \in \mathcal{A}(\mathbb{T})$) if
	
	\begin{align} \label{eq:Ap-condition}
		[p]_{\mathcal{A}} \coloneqq \sup_{J \subset \mathbb{T}} \frac{1}{m(J)} \|\indicator_J\|_{p(\cdot)} \|\indicator_J\|_{p'(\cdot)} < \infty,
	\end{align}
	where the supremum is taken over all subintervals (arcs) $J \subset \mathbb{T} $, and $\indicator_J$ is the indicator function of $J$. This condition is obviously stable under duality, meaning $[p]_{\mathcal{A}} = [p']_{\mathcal{A}}$.
	Recently, in \cite{ADK}, the authors proved that in the case $1<p_{-}\leq p_{+}\leq \infty$, the condition is a necessary and sufficient condition for the boundedness of the Hardy-Littlewood maximal operator in $L^{p(\cdot)}(\mathbb{T})$. Note that if $p_{-}=1$, then the Hardy-Littlewood maximal operator in $L^{p(\cdot)}(\mathbb{T})$ is not bounded (see \cite{CUF}, \cite{DHHR}).
	
	\begin{thm} (\cite{ADK})\label{ADK}
		Let $1<p_{-}\leq p_{+}\leq \infty$. Then the following assertions are equivalent.
		\begin{enumerate}
			\item[(1)] $p(\cdot)\in \mathcal{A}(\mathbb{T})$;
			\item[(2)] The Hardy-Littlewood maximal operator is bounded in $L^{p(\cdot)}(\mathbb{T})$.
		\end{enumerate}
	\end{thm}
	
	To illustrate this result, we will need the following classes of variable exponents. For an exponent $p(\cdot)\in\mathcal{P}(\mathbb{T})$, we say that $p(\cdot) \in \mathcal{P}_{log}(\mathbb{T})$ (logarithmic Dini–Lipschitz condition) if there exists a constant $C$ such that
	$$
	|p(x)-p(y)|\leq C/(-\log |x-y|) \,\,\,\,\mbox{for all} \,\,\,x,y\in \mathbb{T},\,\,\,|x-y|<1/2.
	$$
	
	Note that if $p(\cdot) \in \mathcal{P}_{log}(\mathbb{T})$, then $p(\cdot) \in \mathcal{A}(\mathbb{T})$ (see \cite{CUF}, \cite{DHHR}).
	Under the hypotheses $p_{-}\geq1$ and $p(\cdot) \in \mathcal{P}_{log}(\mathbb{T} )$, Sharapudinov \cite{SH} proved that some families of convolution operators are uniformly bounded in $L^{p(\cdot)}(\mathbb{T})$. This covers in particular a large class of classical operators, including the operators of Fej\'{e}r, de la Vall\'{e}e Poussin, Abel, Steklov and many others.
	
	Lerner \cite{L} constructed an example of a discontinuous variable exponent such that the Hardy-Littlewood maximal operator is bounded on $L^{p(\cdot)}(\mathbb{T})$. Kapanadze and Kopaliani \cite{KK} developed further Lerner’s ideas. In \cite{KK}, the authors considered the following class of variable exponents. Recall that a function $f\in L^{1}(\mathbb{T})$ belongs to the space $VMO^{1/|\log|}$ if
	$$
	\gamma(f,r)=o(1/|\log r|)\,\,\,\,\mbox{as}\,\,\,\,r\rightarrow0,
	$$
	where
	$$
	\gamma(f,r)=\sup_{m(I)\leq r} \frac{1}{m(I)}\int_{I}|f(t)-f_{I}|dm(t).
	$$
	Here $f_{I}$ is the integral average of $f$ on the arc $I$, and the supremum is taken over all arcs $I\subset \mathbb{T}$.
	
	Note that $VMO^{1/|\log|}$ contains discontinuous functions. It is shown in \cite{KK} that if $p\in VMO^{1/|\log|}, \,1<p_{-}\leq p_{+}<\infty$, then the
	Hardy-Littlewood maximal operator is bounded on the variable Lebesgue
	space $L^{p(\cdot)}(\mathbb{T})$.
	
	We say that $w$ is a weight function on $\mathbb{T}$ if $w$ is measurable and $0<w(x)<\infty$ a.e.
	
	\begin{defn}
		We say that $w$ satisfies the $A_{p}(\mathbb{T})$ condition, $1<p<\infty$, if
		$$
		[w]_{A_{p}}=\sup_{I,I\subset\mathbb{T}}\left(\frac{1}{m(I)}\int_{I}w(x)dm(x)\right)\left(\frac{1}{m(I)}\int_{I}w(x)^{-1/(p-1)}dm(x)\right)^{p-1}<\infty.
		$$
		
		We say that $w$ satisfies the $A_{1}(\mathbb{T})$ condition if
		$$
		[w]_{A_{1}}=\sup_{I,I\subset\mathbb{T}}\frac{1}{m(I)}\int_{I}w(x)dm(x) (\mbox{essinf}_{I}w(x))^{-1}<\infty.
		$$
	\end{defn}
	
	Hereafter, we will denote by $\mathcal{S}$ a family of pairs of non-negative,
	measurable functions. Given $p$, $1\leq p<\infty$, if for some $w\in A_{p}(\mathbb{T})$ we write
	$$
	\int_{\mathbb{T}}f(x)^{p}w(x)dm(x) \leq C\int_{\mathbb{T}}g(x)^{p}w(x)dm(x),\,\,\,\,\,(f,g)\in\mathcal{S},
	$$
	then we mean that this inequality holds for all pairs $(f,g)\in\mathcal{S}$ such that the left-hand
	side is finite, and that the constant $C$ may depend on $p$ and $[w]_{A_{p}}$.
	If we write
	$$
	\|f\|_{p(\cdot)}\leq C_{p(\cdot)}\|g\|_{p(\cdot)},\,\,\,\,\,(f,g)\in\mathcal{S},
	$$
	then we mean that this inequality holds for all pairs $(f,g)\in\mathcal{S}$ such that the left-hand side is
	finite and the constant may depend on $p(\cdot)$.
	
	Using this convention, we can state the Rubio
	de Francia extrapolation theorem in the following manner (see \cite{CUF}, Theorem 5.28).
	
	\begin{thm}\label{extrapolation}
		Suppose that for some $p_{0}\geq1$ the family $\mathcal{S}$ is such that for all $w\in A_{1}(\mathbb{T})$
		$$
		\int_{\mathbb{T}}f(x)^{p_{0}}w(x)dm(x)\leq C\int_{\mathbb{T}}g(x)^{p_{0}}w(x)dm(x),\,\,\,\,\,(f,g)\in\mathcal{S}.
		$$
		If, for the exponent $p(\cdot)$, we have $p_{0}\leq p_{-}\leq p_{+}<\infty$ and the maximal operator M is bounded on $L^{(p(\cdot)/p_{0})'}(\mathbb{T})$, then
		$$
		\|f\|_{p(\cdot)}\leq C_{p(\cdot)}\|g\|_{p(\cdot)},\,\,\,\,\,(f,g)\in\mathcal{S}.
		$$
	\end{thm}
	
	The boundedness criterion for Poisson integrals and the Hilbert transform in the $L^{p(\cdot)}(\mathbb{T})$ space may be obtained using the Rubio de Francia extrapolation theorem.
	\begin{thm} \label{PHP}
		Let $1 < p_{-} \leq p_{+} < \infty$. The following statements are equivalent:
		\begin{enumerate}
			\item[(1)] $p(\cdot)\in \mathcal{A}(\mathbb{T})$;
			\item[(2)] $\|P_{r}(f)\|_{p(\cdot)}\leq C \|f\|_{p(\cdot)}$, where the constant $C$ is independent of $r, f.$
			\item[(3)] The circular Hilbert transform is bounded in $L^{p(\cdot)}(\mathbb{T})$.
			\item[(4)] The Riesz projection $P_{+}$ is bounded in $L^{p(\cdot)}(\mathbb{T})$.
		\end{enumerate}
	\end{thm}
	
	Proof. By the theorem of Muckenhoupt, $P_{r}$ is uniformly bounded on the weighted space $L^{p}(\mathbb{T}, w)$, $1<\infty$, if and only if $w$ is in the Muckenhoupt class $A_{p}$ (see \cite{Muc}). Combining the Rubio de Francia extrapolation Theorem \ref{extrapolation} with Theorem \ref{ADK}, we obtain the implication $(1)\Rightarrow(2)$. The implication $(2)\Rightarrow(1)$ may be obtained in a similar way to that for the weighted space $L^{p}(\mathbb{T}, w)$.
	(Note that $\sup_{0\leq r <1}|P_{r}(f,x)|$ is equivalent to the Hardy-Littlewood maximal function for a nonnegative function $f$).
	
	By the theorem of Hunt-Muckenhoupt-Wheeden, the circular Hilbert transform is bounded on the weighted space $L^{p}(\mathbb{T}, w)$, $1<p<\infty$, if and only if $w$ is in the Muckenhoupt class $A_{p}$ (see \cite{Gar}). Combining the Rubio de Francia extrapolation Theorem \ref{extrapolation} with Theorem \ref{ADK}, we obtain the implication $(1)\Rightarrow(3)$. The implication $(3)\Rightarrow (1)$ may be obtained in a similar way to that for the weighted space $L^{p}(\mathbb{T}, w)$.
	
	Note that $P_{+}=\frac{1}{2}(I+iH+P_{0})$, where $I$ is the identity operator and $P_{0}f=\widehat{f}(0)$ (for the definition of $\widehat{f}(0)$, see below), and consequently $(1)\Leftrightarrow(4)$ follows from $(1)\Leftrightarrow(3)$.  $\Box$
	
	By the theorem of Muckenhoupt, $P_{r}$ is uniformly bounded on the weighted space $L^{1}(\mathbb{T})$ if and only if $w$ is in the Muckenhoupt class $A_{1}$ (see \cite{Gar}). The boundedness criterion for Poisson integrals in the $L^{p(\cdot)}(\mathbb{T})$ space in the case $1 \leq p_{-} \leq p_{+} < \infty$ may be obtained analogously using the Rubio de Francia extrapolation theorem.
	
	\begin{thm} \label{P1}
		Let $1 \leq p_{-} \leq p_{+} < \infty$. The following statements are equivalent:
		\begin{enumerate}
			\item[(1)] $p(\cdot)\in \mathcal{A}(\mathbb{T})$;
			\item[(2)] $\|P_{r}(f)\|_{p(\cdot)}\leq C \|f\|_{p(\cdot)}$, where the constant $C$ is independent of $r, f.$
		\end{enumerate}
	\end{thm}
	
	\begin{cor}\label{limit}
		Let $1 \leq p_{-} \leq p_{+} < \infty$ and $p(\cdot)\in \mathcal{A}(\mathbb{T})$.
		Then $\|f-P_{r}f\|_{p(\cdot)}\rightarrow 0$ as $r\rightarrow1^{-}.$
	\end{cor}
	
	\section{Variable exponent Hardy space in the unit disk}
	
	The space $H^{p(\cdot)}(\mathbb{D})$ is the collection of functions $F(z)=F(re^{it})$ that are analytic in the unit disk and have finite norm
	$$
	\|F\|_{H^{p(\cdot)}(\mathbb{D})}=\sup _{0<r<1} \|F_{r}\|_{p(\cdot)},
	$$
	where $F_{r}:\mathbb{T}\rightarrow \mathbb{C}$ is defined as $F_{r}(t)=F(re^{it}).$
	
	For $f\in L^{1}(\mathbb{T})$, let
	$$
	\widehat{f}(n)=\frac{1}{2\pi}\int_{0}^{2\pi}f(e^{it})e^{-int}dt,\,\,\,n\in \mathbb{Z}.
	$$
	
	We will consider the Hardy space $H^{p(\cdot)}(\mathbb{T})$ built upon the space $L^{p(\cdot)}(\mathbb{T})$, which is defined by
	$$
	H^{p(\cdot)}(\mathbb{T})=\left\lbrace f\in L^{p(\cdot)}(\mathbb{T}:\widehat{f}(n)=0\,\,\,\mbox{for all}\,\,n<0 \right\rbrace .
	$$
	
	It is clear that $H^{p(\cdot)}(\mathbb{T})$ is a closed subspace of $L^{p(\cdot)}(\mathbb{T}).$
	
	A function of the form
	$$
	\sum_{k=0}^{n} a_{k}t^{k},\,\,t\in \mathbb{T},\,\,\,a_{0},...,a_{n}\in \mathbb{C}
	$$
	is said to be an analytic polynomial on $\mathbb{T}$. The set of all analytic polynomials is denoted by $\varPi_{\mathbb{A}}$. It is well known that the set $\varPi_{\mathbb{A}}$ is dense in $H^{p}(\mathbb{T})$ when $p(x)=p=const,\,1\leq p<\infty$. Using Theorem 1.2 from \cite{Kar}, we obtain that
	if $1< p_{-}\leq p_{+}<\infty$ and $p(\cdot) \in \mathcal{A}(\mathbb{T})$, then $\varPi_{\mathbb{A}}$ is dense in $ H^{p(\cdot)}(\mathbb{T})$.
	
	As in the classical case, when $1\leq p_{-}\leq p_{+}<\infty$ and $p(\cdot) \in \mathcal{A}(\mathbb{T})$, we may identify a function $f\in H^{p(\cdot)}(\mathbb D)$ with its boundary values $bf\in H^{p(\cdot)}(\mathbb T)$, and hence also the spaces $H^{p(\cdot)}(\mathbb D)$
	and $H^{p(\cdot)}(\mathbb T)$. In particular, we have
	
	\begin{thm}\label{boundary-value-theorem}
		Let $1\leq p_{-}\leq p_{+}<\infty, \qquad p(\cdot)\in\mathcal A(\mathbb T)$, and let $F\in H^{p(\cdot)}(\mathbb D)$.
		Then there exists a function $\Phi\in H^{p(\cdot)}(\mathbb T)$ such that
		\begin{equation}\label{ae-boundary-limit}
			\lim_{r\to1^-}F(re^{ix})=\Phi(x)
		\end{equation}
		for almost every $x\in[-\pi,\pi]$.
		Moreover:
		
		\begin{enumerate}
			\item[(1)] For $z=re^{ix}\in\mathbb D$,
			\begin{equation}\label{Poisson-representation-boundary}
				F(re^{ix})
				=
				\frac{1}{2\pi}
				\int_{-\pi}^{\pi}
				\Phi(x-t)P_r(t)\,dt,
			\end{equation}
			
			\item[(2)]
			\begin{equation}\label{Lp-radial-convergence}
				\lim_{r\to1^-}
				\|F(re^{i\cdot})-\Phi\|_{p(\cdot)}
				=0.
			\end{equation}
			
			\item[(3)] There exists a constant
			$C_{p(\cdot)}\geq1$, depending only on $p(\cdot)$, such that
			\begin{equation}\label{Hardy-boundary-equivalence}
				\|\Phi\|_{p(\cdot)}
				\leq
				\|F\|_{H^{p(\cdot)}(\mathbb D)}
				\leq
				C_{p(\cdot)}
				\|\Phi\|_{p(\cdot)}.
			\end{equation}
		\end{enumerate}
	\end{thm}
	
	Proof. Since $H^{p(\cdot)}(\mathbb D)\subset H^1(\mathbb D)$, we have $F\in H^1(\mathbb D)$ and, by the classical Fatou theorem for $H^1$, there exists
	$\Phi\in L^1(\mathbb T)$ such that
	\begin{equation}\label{Fatou-limit-proof}
		F(re^{ix})\longrightarrow\Phi(x)
	\end{equation}
	for almost every $x\in[-\pi,\pi]$ as $r\to1^-$. Furthermore, $\Phi$ has no negative Fourier coefficients:
	\begin{equation}\label{negative-coeff-zero}
		\widehat{\Phi}(n)=0,
		\qquad n<0.
	\end{equation}
	
	Let $\lambda>\|F\|_{H^{p(\cdot)}(\mathbb D)}$.
	Using the almost everywhere convergence
	\eqref{Fatou-limit-proof} and Fatou's lemma, we obtain
	\[
	\rho_{p(\cdot)}
	\left(\frac{\Phi}{\lambda}\right)=\int_{\mathbb T}\left|\frac{\Phi(x)}{\lambda}\right|^{p(x)}dm(x)
	\leq \liminf_{r\to1^-}\int_{\mathbb T}\left|\frac{F_r(x)}{\lambda}\right|^{p(x)}dm(x)
	\leq1.
	\]
	
	Together with \eqref{negative-coeff-zero}, this shows that $ \Phi\in H^{p(\cdot)}(\mathbb T)$ and
	
	\begin{equation}\label{boundary-norm-lower}
		\|\Phi\|_{p(\cdot)}
		\leq \|F\|_{H^{p(\cdot)}(\mathbb D)}.
	\end{equation}
	
	Since $F\in H^1(\mathbb D)$ and $\Phi$ is its boundary function,
	the classical boundary representation theorem for $H^1$ gives
	\begin{equation}\label{classical-Poisson-representation}
		F(re^{ix})
		=
		P_r*\Phi(x).
	\end{equation}
	This proves assertion $(1)$.
	
	Applying Corollary \ref{limit} and using
	\eqref{classical-Poisson-representation}, we obtain
	$$
	\|F(re^{i\cdot})-\Phi\|_{p(\cdot)}=\|P_r*\Phi-\Phi\|_{p(\cdot)} \longrightarrow0.
	$$
	This proves assertion $(2)$.
	
	Since
	\[
	F_r=P_r*\Phi,
	\]
	by the uniform boundedness of the
	Poisson operators on $L^{p(\cdot)}(\mathbb T)$, we obtain
	\[
	\|F_r\|_{p(\cdot)}
	\leq
	C_{p(\cdot)}
	\|\Phi\|_{p(\cdot)}
	\]
	for every $0<r<1$. Together with \eqref{boundary-norm-lower}, this proves
	\eqref{Hardy-boundary-equivalence} and assertion $(3)$.
	
	Notice also that assertion $(2)$ implies $ \|F_r\|_{p(\cdot)} \longrightarrow \|\Phi\|_{p(\cdot)}$ as $r\to1^-$.  $\Box$
	
	Note that for a variable exponent the norms $\|F_r\|_{p(\cdot)}$ need not be
	monotone in $r$, and the equality
	\[
	\sup_{0<r<1}\|F_r\|_{p(\cdot)}
	=
	\|\Phi\|_{p(\cdot)}
	\]
	is in general false for nonconstant variable exponents.
	
	For example, let
	\[
	p(t)=2+\delta\cos t,
	\qquad 0<\delta<1.
	\]
	Then
	\[
	1<2-\delta\leq p(t)\leq2+\delta<\infty,
	\]
	and $p$ is smooth; in particular,
	$p(\cdot)\in\mathcal A(\mathbb T)$.
	
	For $\varepsilon>0$, define $F_\varepsilon(z)=1-\varepsilon z$.
	Its boundary function is $\Phi_\varepsilon(t)=1-\varepsilon e^{it}$.
	
	Consider its modular at $\lambda=1$:
	\[
	I(\varepsilon)
	=
	\rho_{p(\cdot)}(\Phi_\varepsilon)
	=
	\frac1{2\pi}
	\int_{-\pi}^{\pi}
	|1-\varepsilon e^{it}|^{2+\delta\cos t}\,dt.
	\]
	Clearly, $I(0)=1$. Moreover,
	\[
	|1-\varepsilon e^{it}|^{p(t)}
	=
	\left(
	1-2\varepsilon\cos t+\varepsilon^2
	\right)^{p(t)/2},
	\]
	and therefore
	\[
	I'(0)=-\frac{1}{2\pi}\int_{-\pi}^{\pi}p(t)\cos t\,dt=-\frac{1}{2\pi}\int_{-\pi}^{\pi}
	(2+\delta\cos t)\cos t\,dt=-\frac{\delta}{2}<0.
	\]
	
	Hence, for all sufficiently small $\varepsilon>0$, $I(\varepsilon)<1.$
	By continuity of the modular with respect to the scaling parameter,
	this implies $ \|\Phi_\varepsilon\|_{p(\cdot)}<1.$
	
	On the other hand, $ F_\varepsilon(0)=1.$
	
	Moreover, $ F_{\varepsilon}(re^{it})\longrightarrow1 $
	uniformly as $r\to0^+$, and hence
	\[
	\|F_\varepsilon(re^{i\cdot})\|_{p(\cdot)}
	\longrightarrow
	\|1\|_{p(\cdot)}
	=1.
	\]
	Therefore
	\[
	\|F_\varepsilon\|_{H^{p(\cdot)}(\mathbb D)}
	=
	\sup_{0<r<1}
	\|F_\varepsilon(re^{i\cdot})\|_{p(\cdot)}
	\geq1
	>
	\|\Phi_\varepsilon\|_{p(\cdot)}.
	\]
	Thus exact norm equality cannot in general hold for nonconstant
	$p(\cdot)$.
	
	\section {Proof of Theorem \ref{main Theorem1}}
	
	\medskip
	\noindent
	\textbf{Proof of $(i)\Rightarrow(ii)$.}
	
	For $f\in H^{p(\cdot)}(\mathbb T)$, set
	\[
	S_Nf
	=
	\sum_{n=0}^{N-1}
	\langle f,\phi_n\rangle\phi_n.
	\]
	
	Note that
	\begin{equation}\label{MT-kernel}
		\sum_{n=0}^{N-1}
		\phi_n(z)\overline{\phi_n(w)}
		=
		\frac{
			1-B_N(z)\overline{B_N(w)}
		}
		{1-z\overline w}.
	\end{equation}
	
	Since $p_->1$, we have
	$H^{p(\cdot)}(\mathbb T)\subset H^1(\mathbb T)$.
	Using \eqref{MT-kernel} and the Cauchy formula, for $z\in\mathbb D$
	we obtain
	\begin{align*}
		S_Nf(z)
		&=
		\int_{\mathbb T}
		f(\zeta)
		\frac{
			1-B_N(z)\overline{B_N(\zeta)}
		}
		{1-z\overline\zeta}
		\,dm(\zeta)
		\\
		&=
		f(z)
		-
		B_N(z)
		\int_{\mathbb T}
		\frac{
			\overline{B_N(\zeta)}f(\zeta)
		}
		{1-z\overline\zeta}
		\,dm(\zeta).
	\end{align*}
	Hence, taking boundary values, we obtain
	\begin{equation}\label{SN-Riesz}
		S_Nf
		=
		f-B_NP_+(\overline{B_N}f).
	\end{equation}
	
	Since $|B_N|=1$ almost everywhere on $\mathbb T$, multiplication by
	$B_N$ or by $\overline{B_N}$ is an isometry of
	$L^{p(\cdot)}(\mathbb T)$. Therefore, by Theorem
	\ref{PHP} and \eqref{SN-Riesz},
	\begin{align}
		\|S_Nf\|_{p(\cdot)}
		&\leq
		\|f\|_{p(\cdot)}
		+
		\|P_+(\overline{B_N}f)\|_{p(\cdot)}
		\nonumber\\
		&\leq
		(1+C_{p(\cdot)})\|f\|_{p(\cdot)}.
		\label{SN-uniform}
	\end{align}
	
	Let
	$ \mathcal D =\operatorname{span} \{k_a:a\in\mathbb D\},$
	where $ k_a(z)=\frac1{1-\overline a z}. $
	We claim that $\mathcal D$ is dense in
	$H^{p(\cdot)}(\mathbb T)$. Indeed, since $p(\cdot) \in \mathcal{A}(\mathbb{T})$, $\varPi_{\mathbb{A}}$ is dense in $ H^{p(\cdot)}(\mathbb{T})$. The systems $\mathcal D$ and $\varPi_{\mathbb{A}}$ are dense in $H^{p_{+}}(\mathbb{T})$. Using the fact that $H^{p_{+}}(\mathbb{T}) \subset H^{p(\cdot)}(\mathbb{T})$, we obtain that
	$\mathcal D$ is dense in the space $ H^{p(\cdot)}(\mathbb{T})$.
	
	To prove the implication $(i)\Rightarrow(ii)$, it suffices to show that
	$$
	\|S_Nk_a-k_a\|_{p(\cdot)}
	\longrightarrow0
	\qquad(a\in\mathbb D).
	$$
	
	By the reproducing property, we have $\langle k_a,\phi_n\rangle=\overline{\phi_n(a)}$ and, by \eqref{MT-kernel},
	$$
	S_Nk_a(z)=\sum_{n=0}^{N-1}\overline{\phi_n(a)}\phi_n(z)=\frac{
		1-B_N(z)\overline{B_N(a)}}{1-z\overline a}.
	$$
	Consequently, $ k_a(z)-S_Nk_a(z)=\overline{B_N(a)}\,B_N(z)k_a(z).$ Since $|B_N|=1$ on $\mathbb T$ and $B_N(a)\longrightarrow0$
	$ (a\in\mathbb D)$, we obtain
	$$
	\|k_a-S_Nk_a\|_{p(\cdot)}=|B_N(a)|\,\|k_a\|_{p(\cdot)}\rightarrow 0, \,\,\,\mbox{as}\,\,\, N\rightarrow\infty.
	$$
	
	\medskip
	\noindent
	\textbf{Proof of $(ii)\Rightarrow(i)$.}
	
	Assume now that $(\phi_n)_{n\geq0}$ is a Schauder basis of
	$H^{p(\cdot)}(\mathbb T)$. We will prove that the projection operator $P_{+}$ is bounded on $H^{p(\cdot)}(\mathbb T)$. From Theorem \ref{PHP}, we obtain $p(\cdot)\in \mathcal{A}(\mathbb{T})$.
	
	The uniform boundedness principle gives
	\begin{equation}\label{basis-uniform}
		K=
		\sup_{N\geq1}
		\|S_N\|_{H^{p(\cdot)}\to H^{p(\cdot)}}
		<\infty,
	\end{equation}
	where $S_Nf$ are the partial sum operators of $f\in H^{p(\cdot)}(\mathbb T).$
	
	For $N\geq1$, define $ X_N=\overline{B_N}H^{p(\cdot)}(\mathbb T)$.
	
	Observe that $X_N\subset X_{N+1}$. Indeed, since $B_{N+1}=B_Nb_{a_{N+1}}$, if $g=\overline{B_N}f$, then on $\mathbb T$
	$ g=\overline{B_{N+1}} \bigl(b_{a_{N+1}}f\bigr)$, and $b_{a_{N+1}}f\in H^{p(\cdot)}(\mathbb T)$.
	
	Note that for every $N$, $H^{p(\cdot)}(\mathbb T)\subset X_N$, because $h=\overline{B_N}(B_Nh),$ $ h\in H^{p(\cdot)}$.
	
	We now prove that the set $ \cup_{N=1}^{\infty}X_N$ is dense in $L^{p(\cdot)}(\mathbb T)$. Suppose, to the contrary, that this is false.
	By the Hahn--Banach theorem, there exists a non-zero
	$g\in L^{p'(\cdot)}(\mathbb T)$ such that
	\begin{equation}\label{annihilator}
		\int_{\mathbb T}
		h\,\overline g\,dm
		=
		0
		\qquad
		\text{for every }
		h\in\bigcup_{N=1}^{\infty}X_N.
	\end{equation}
	
	Note that \eqref{annihilator} implies
	\[
	\int_{\mathbb T}
	z^k\overline g\,dm=0,
	\qquad k\geq0.
	\]
	Equivalently,
	\begin{equation}\label{g-negative}
		\widehat g(k)=0,
		\qquad k\geq0.
	\end{equation}
	
	Define
	$G(\zeta)=\overline\zeta\,\overline{g(\zeta)},\qquad \zeta\in\mathbb T.$
	By \eqref{g-negative}, $G$ has no negative Fourier coefficients.
	Moreover, $ |G|=|g|,$ and $ G\in H^{p'(\cdot)}(\mathbb T).$
	
	Fix $N$. Since
	$\overline{B_N}h\in X_N$ for every
	$h\in H^{p(\cdot)}(\mathbb T)$, \eqref{annihilator} gives
	\[
	\int_{\mathbb T}
	\overline{B_N}h\,\overline g\,dm=\int_{\mathbb T} h\,\overline{B_Ng}\,dm=0
	\]
	and consequently $\widehat{B_Ng}(k)=0,\,\, k\geq0$.
	The function
	\[
	G_N(\zeta)
	:=
	\overline\zeta\,
	\overline{B_N(\zeta)g(\zeta)}
	=
	\overline{B_N(\zeta)}G(\zeta)
	\]
	has no negative Fourier coefficients. Since
	$|G_N|=|G|$, we have $ G_N\in H^{p'(\cdot)}(\mathbb T).$
	
	Note that on $\mathbb T$, $G=B_NG_N $ and both $G$ and $G_N$ belong to $H^1(\mathbb T)$. Hence their analytic
	extensions satisfy
	\begin{equation}\label{divisibility}
		G(z)=B_N(z)G_N(z),
		\qquad z\in\mathbb D.
	\end{equation}
	
	Since \eqref{divisibility} holds for every $N$, the function $G$
	vanishes at every point $a_n$, with the corresponding multiplicity.
	If $G$ were not identically zero, then the zero sequence of the
	non-zero $H^1$ function $G$ would satisfy the classical Blaschke
	condition
	\[
	\sum_{n=1}^{\infty}(1-|a_n|)<\infty;
	\]
	see, for example, \cite{Du}. This contradicts
	\eqref{Blaschke condition}. Hence $G\equiv0$, and therefore
	$g=0$, contradicting the choice of $g$. Consequently, the set $ \cup_{N=1}^{\infty}X_N$ is dense in $L^{p(\cdot)}(\mathbb T)$.
	
	If $g\in X_N$, write $ g=\overline{B_N}f,\, \mbox{where}\, f\in H^{p(\cdot)}(\mathbb T)$.
	
	By \eqref{MT-kernel}, for $z\in\mathbb D$,
	\begin{align*}
		S_Nf(z)
		&=
		f(z)
		-
		B_N(z)
		\int_{\mathbb T}
		\frac{g(\zeta)}
		{1-z\overline\zeta}
		\,dm(\zeta).
	\end{align*}
	Hence
	\begin{equation}\label{Cauchy-positive}
		\frac{f(z)-S_Nf(z)}{B_N(z)}
		=
		\int_{\mathbb T}
		\frac{g(\zeta)}
		{1-z\overline\zeta}
		\,dm(\zeta)
		=
		\sum_{k=0}^{\infty}
		\widehat g(k)z^k.
	\end{equation}
	
	Consequently, the boundary function $T_Ng= \overline{B_N}(f-S_Nf)$ has Fourier coefficients
	\begin{equation}\label{TN-Fourier}
		\widehat{T_Ng}(k)
		=
		\begin{cases}
			\widehat g(k),& k\geq0,\\
			0,& k<0.
		\end{cases}
	\end{equation}
	
	Moreover, using \eqref{basis-uniform},
	\begin{align}
		\|T_Ng\|_{p(\cdot)}
		&=
		\|f-S_Nf\|_{p(\cdot)}
		\nonumber\\
		&\leq
		(1+K)\|f\|_{p(\cdot)}
		=
		(1+K)\|g\|_{p(\cdot)}.
		\label{TN-bound}
	\end{align}
	
	Put
	\[
	\mathcal X
	=
	\bigcup_{N=1}^{\infty}X_N.
	\]
	Because the spaces $X_N$ are increasing, $\mathcal X$ is a linear
	subspace. Define a linear operator
	\[
	T:\mathcal X\longrightarrow L^{p(\cdot)}(\mathbb T)
	\]
	such that
	\[
	\widehat{Tg}(k)
	=
	\begin{cases}
		\widehat g(k),&k\geq0,\\
		0,&k<0,
	\end{cases}
	\]
	and, by \eqref{TN-bound},
	\[
	\|Tg\|_{p(\cdot)}
	\leq
	(1+K)\|g\|_{p(\cdot)}.
	\]
	Since $\mathcal X$ is dense in $L^{p(\cdot)}(\mathbb T)$, $T$ extends uniquely to a bounded operator on
	all of $L^{p(\cdot)}(\mathbb T)$.
	
	We claim that this extension is the Riesz projection $P_+$. Indeed,
	let $g_j\in\mathcal X$ and $ g_j\longrightarrow g \,\,\text{in }L^{p(\cdot)}(\mathbb T)$.
	
	Since  $ L^{p(\cdot)}(\mathbb T)\hookrightarrow L^1(\mathbb T)$ we have 
	$ \widehat{g_j}(k)\longrightarrow\widehat g(k) $
	for every $k\in\mathbb Z$. Also $Tg_j\longrightarrow Tg$ in $L^{p(\cdot)}$ and therefore in $L^1$. Consequently 
	\[
	\widehat{Tg}(k)
	=
	\begin{cases}
		\widehat g(k),&k\geq0,\\
		0,&k<0.
	\end{cases}
	\]
	
	Thus $T=P_+$, and therefore
	\begin{equation}\label{Pplus-necessary}
		P_+:
		L^{p(\cdot)}(\mathbb T)
		\longrightarrow
		L^{p(\cdot)}(\mathbb T)
	\end{equation}
	is bounded.
	
	It follows that the circular Hilbert transform is bounded on
	$L^{p(\cdot)}(\mathbb T)$, and consequently
	$p(\cdot)\in\mathcal A(\mathbb T)$.    $\Box$
	
	\section {Proof of Theorem \ref{equivalence}}
	
	Since $ B_n+\overline{a_{n+1}}B_{n+1}=\sqrt{1-|a_{n+1}|^2}\,\phi_n$,
	we have
	\begin{equation}\label{psi-B-relation}
		\psi_n
		=
		B_n+\overline{a_{n+1}}B_{n+1},
	\end{equation}
	where $ \psi_n=\sqrt{1-|a_{n+1}|^2}\,\phi_n$.
	
	Note that $\sqrt{1-r^2}\leq \sqrt{1-|a_{n+1}|^2}\leq1$, and consequently the system $(\psi_n)$ is a Schauder basis if and only if
	$(\phi_n)$ is a Schauder basis. Thus it remains to prove that
	$(B_n)$ is a Schauder basis if and only if $(\psi_n)$ is a
	Schauder basis.
	
	Assume that $(B_n)_{n\geq0}$ is a Schauder basis. For $f\in H^{p(\cdot)}(\mathbb T)$, there is a unique sequence
	$(c_n)_{n\geq0}$ such that
	\begin{equation}\label{B-expansion}
		f=\sum_{n=0}^{\infty}c_nB_n
	\end{equation}
	with convergence in $H^{p(\cdot)}(\mathbb T)$. Since $|B_n(\zeta)|=1,
	\, (\zeta\in\mathbb T)$, we have $ \|B_n\|_{p(\cdot)}=1,$
	and since the series \eqref{B-expansion} converges, we obtain $ c_nB_n\longrightarrow0
	$ in $H^{p(\cdot)}$. Consequently, we have $c_n\longrightarrow0$.
	Set $ \lambda_n:=\overline{a_n},
	\,(n\geq1)$ and define $(d_n)_{n\geq0}$ recursively by
	\begin{equation}\label{d-recurrence}
		d_0=c_0,
		\qquad
		d_n=c_n-\lambda_nd_{n-1},
		\quad n\geq1.
	\end{equation}
	Equivalently,
	\begin{equation}\label{c-d-relation}
		c_0=d_0,
		\qquad
		c_n=d_n+\lambda_nd_{n-1},
		\quad n\geq1.
	\end{equation}
	
	Iterating \eqref{d-recurrence} and using $|\lambda_n|\leq r$, we
	obtain
	\begin{equation}\label{dn-estimate}
		|d_n|
		\leq
		|c_n|
		+r|c_{n-1}|
		+r^2|c_{n-2}|
		+\cdots
		+r^n|c_0|.
	\end{equation}
	We claim that $ d_n\longrightarrow0. $
	Indeed, let $\varepsilon>0$. Choose $N_0$ such
	that
	\[
	|c_j|<\varepsilon(1-r)
	\qquad(j\geq N_0).
	\]
	For $n\geq N_0$, split the right-hand side of
	\eqref{dn-estimate} into
	\[
	\sum_{j=N_0}^{n}
	r^{\,n-j}|c_j|
	+
	\sum_{j=0}^{N_0-1}
	r^{\,n-j}|c_j|.
	\]
	For the first sum,
	\[
	\sum_{j=N_0}^{n}
	r^{\,n-j}|c_j|
	\leq
	\varepsilon(1-r)
	\sum_{k=0}^{\infty}r^k
	=
	\varepsilon.
	\]
	The second sum tends to zero as $n\to\infty$, since it contains only
	finitely many fixed numbers multiplied by powers of $r$. Consequently, $ d_n\longrightarrow0.$
	
	Using \eqref{psi-B-relation} and \eqref{c-d-relation}, we obtain,
	for every $N\geq0$,
	\begin{align}
		\sum_{n=0}^{N}d_n\psi_n
		&=
		\sum_{n=0}^{N}
		d_n
		\left(
		B_n+\lambda_{n+1}B_{n+1}
		\right)
		\nonumber\\
		&=
		d_0B_0
		+
		\sum_{n=1}^{N}
		\left(d_n+\lambda_nd_{n-1}\right)B_n
		+
		\lambda_{N+1}d_NB_{N+1}
		\nonumber\\
		&=
		\sum_{n=0}^{N}c_nB_n
		+
		\lambda_{N+1}d_NB_{N+1}.
		\label{partial-sum-relation}
	\end{align}
	 We obtain
	\[
	\|\lambda_{N+1}d_NB_{N+1}\|_{p(\cdot)}
	=
	|\lambda_{N+1}||d_N|
	\leq
	r|d_N|
	\longrightarrow0.
	\]
	Passing to the limit in \eqref{partial-sum-relation}, we obtain
	\begin{equation}\label{psi-expansion}
		f=\sum_{n=0}^{\infty}d_n\psi_n.
	\end{equation}
	Thus every $f\in H^{p(\cdot)}(\mathbb T)$ has an expansion with
	respect to $(\psi_n)$. The uniqueness of the representation follows easily from the embedding  $H^{p(\cdot)}(\mathbb{T})\subset H^{p_{-}}(\mathbb{T})$ considering that the system $(\phi_n)_{n\geq0}$  is a basis for the space $H^{p_{-}}(\mathbb{T})$ (see \cite{CP}).
	
	Assume that $(\psi_n)_{n\geq0}$ is a Schauder basis.
	Let $f\in H^{p(\cdot)}(\mathbb T)$ and write
	\begin{equation}\label{psi-expansion2}
		f=\sum_{n=0}^{\infty}d_n\psi_n.
	\end{equation}
	
	We first show that $d_n\to0$. On $\mathbb T$, since
	$ |\psi_n| \geq 1-|a_{n+1}| \geq 1-r$, we have
	\begin{equation}\label{psi-lower-bound}
		\|\psi_n\|_{p(\cdot)}
		\geq
		(1-r)\|1\|_{p(\cdot)}
		=
		1-r.
	\end{equation}
	Since the series \eqref{psi-expansion2} converges, $ d_n\psi_n\longrightarrow0$ in norm.
	Hence, by \eqref{psi-lower-bound}, $(1-r)|d_n| \leq \|d_n\psi_n\|_{p(\cdot)}\longrightarrow0, $
	and consequently $ d_n\longrightarrow0 $.
	Define
	\begin{equation}\label{define-cn-converse}
		c_0=d_0,
		\qquad
		c_n=d_n+\lambda_nd_{n-1},
		\quad n\geq1.
	\end{equation}
	As in \eqref{partial-sum-relation},
	\[
	\sum_{n=0}^{N}d_n\psi_n
	=
	\sum_{n=0}^{N}c_nB_n
	+
	\lambda_{N+1}d_NB_{N+1}.
	\]
	Therefore
	\begin{equation}\label{B-partial-from-psi}
		\sum_{n=0}^{N}c_nB_n
		=
		\sum_{n=0}^{N}d_n\psi_n
		-
		\lambda_{N+1}d_NB_{N+1}.
	\end{equation}
	Since
	\[
	\|\lambda_{N+1}d_NB_{N+1}\|_{p(\cdot)}
	\leq
	r|d_N|
	\longrightarrow0,
	\]
	we obtain
	\begin{equation}\label{B-expansion-converse}
		f=\sum_{n=0}^{\infty}c_nB_n.
	\end{equation}
	Thus every $f\in H^{p(\cdot)}(\mathbb T)$ has a convergent expansion
	with respect to $(B_n)$. The uniqueness of the representation follows easily from the embedding $H^{p(\cdot)}(\mathbb{T})\subset H^{p_{-}}(\mathbb{T})$ considering that the system $(B_n)_{n\geq0}$ is a basis for the space $H^{p_{-}}(\mathbb{T})$ (see \cite{FMN}).
	 $\Box$

\end{document}